\documentclass[letterpaper, 10 pt, conference]{ieeeconf}  % Comment this line out if you need a4paper

\IEEEoverridecommandlockouts                              % This command is only needed if 
\usepackage{graphics} % for pdf, bitmapped graphics files
\usepackage{epsfig} % for postscript graphics files
\usepackage{mathptmx} % assumes new font selection scheme installed
\usepackage{times} % assumes new font selection scheme installed
\usepackage{amsmath} % assumes amsmath package installed
\usepackage{amssymb}  % assumes amsmath package installed

\usepackage{amsrefs}
\usepackage{xcolor}

\newtheorem{theorem}{Theorem}[section]
\newtheorem{fact}[theorem]{Fact}
\newcommand{\bfact}{\begin{fact}}
	\newcommand{\efact}{\end{fact}}
\newtheorem{proposition}[theorem]{Proposition}

\newtheorem{definition}[theorem]{Definition}
\newcommand{\sgn}{\text{sgn}}
\newtheorem{example}[theorem]{Example}

\newtheorem{remark}[theorem]{Remark}

\def\rr{\mathbb{R}}
\def\nn{\mathbb{N}}

\def\Lin{\text{Lin}}

\def\K{\mathcal{K}}

\def\M{\mathcal{M}}
\def\N{\mathcal{N}}

\def\Span{\text{Span}}
\def\Lin{\text{Lin}}
\def\ds{\displaystyle}
\def\bel{\begin{equation}\label}
\def\eeq{\end{equation}}
\def\fra {\color{blue}}

\def\E{\mathcal{E}}
\usepackage{esint}

\title{\LARGE \bf
Quasi Differential Quotients
}

\author{Francesca Angrisani$^{\dagger}$ and Franco Rampazzo$^{\ddagger}$% <-this % stops a space
\thanks{$^{\dagger}$ Francesca Angrisani is with Department of Mathematics,
        Universita' degli Studi di Padova, 35131 Padova, Italy
        {\tt\small francesca.angrisani@unina.it}}%
\thanks{$^{\ddagger}$Franco Rampazzo is with the Department of Mathematics, Universita' degli Studi di Padova, 35131 Padova, Italy
        {\tt\small rampazzo.math@unipd.it}}%
}

\binoppenalty=\maxdimen
\relpenalty=\maxdimen
\begin{document}

\maketitle
\thispagestyle{empty}
\pagestyle{empty}

%%%%%%%%%%%%%%%%%%%%%%%%%%%%%%%%%%%%%%%%%%%%%%%%%%%%%%%%%%%%%%%%%%%%%%%%%%%%%%%%
\begin{abstract}

	We explore basic properties and some applications of Quasi Differential Quotients ($QDQ$s) and the related  $QDQ$-approximating multi-cones. A $QDQ$, which is a special kind of   H.Sussmann's Approximate Generalized Differential Quotient ($AGDQ$), consists in a notion of generalized differentiation   for  set-valued maps. $QDQ$s have the advantage over $AGDQ$s of allowing a genuine, non-punctured, Open Mapping result, so implying stronger set-separation theorems. They have  already proved  quite useful  in the investigation  of 
	some connections occurring between   infimum gap phenomena  and the  normality of minima. Moreover,  $QDQ$-approximating  multi-cones are fit in optimal control to deduce Maximum Principles  that involve  (set-valued) Lie  brackets of nonsmooth vector fields.   
	
\end{abstract}

%%%%%%%%%%%%%%%%%%%%%%%%%%%%%%%%%%%%%%%%%%%%%%%%%%%%%%%%%%%%%%%%%%%%%%%%%%%%%%%%
\section{INTRODUCTION}
{\fra

}

Extending a result due to Carath\'eodory,  M.W. Botsko \cite{Botsko} gave the following  simple characterization of differentiability: a function $f:\rr^n\to\rr$ is differentiable at $a\in\rr^n$ if and only if the exists $\bar\delta$ and  a map $L:a+B_{\bar\delta}\to (\rr^n)^*$, continuous at $a$ such that
$$ f(x) = f(a) + L(x) (x-a) \quad \forall x\in a+B_{\bar\delta}. $$
On the same line of thought, Sussmann proposed (see \cite{sussW}) two instances  of generalized differentials for set-valued maps, the ``Generalized Differential Quotients ($GDQ$)"  and the ``Approximate Generalized Differential Quotient" ($AGDQ$).  This was part of a program designed to identifying generalized differentiation theories capable of {\it automatically}  providing  a Maximum Principle for minimizers of optimal control problems. This program naturally consists   in  requiring  some properties to be verified by the proposed generalized differentials, including
{\it locality, product rules, chain rule, and open mapping results for the related approximating multi-cones.} 
Recently, a special subclass of $AGDQ$s has been proposed in \cite{PR} in order to establish a general, {\it normality}-based, criterion for the avoidance of infimum-gaps. These particular $AGDQ$s have been called {\it Quasi Differential Quotients} (QDQs).  Let us also mention that $AGDQ$, hence  $QDQ$s, and the related approximating  multi-cones are chart-independent objects, so that  they  can be defined on differential manifolds.The main property of $QDQ$s consists in the validity of a genuine open mapping result, whereas only a ``punctured"  open mapping holds true for general $AGDQ$s. This fact has proved essential in the proof of the quoted criterion. At the same time,
$QDQ$s are going to be utilized for other applications. In particular, they prove to be a quite ductile tool in the program of extending Goh-like conditions for optimal control problems with nonsmooth dynamics. Such conditions involve set-valued Lie brackets, which happen to be $QDQ$s of suitable compositions of flows. 

In the present paper we wish to present the notion of $QDQ$, by proving  some basic properties as well as by recalling some  known facts. We also provide some significant examples and give sufficient or necessary  conditions  for a set to be a $QDQ$. This is intended as the beginning of a more articulate investigation  addressed on one hand  to the relation between  $QDQ$s with other notions of generalized differentiation (see e.g. \cite{Mordukhovich}, \cite{Clarke})  and, on the other hand, to concrete applications.
Here is an outline of the article. In   Section 2, we present the definition of $QDQ$s in Euclidean spaces and on manifolds, and prove basic properties, comprising  locality, linearity, product rules,  and a Chain Rule. Section 3 is devoted to $QDQ$-approximating multi-cones and to the above-mentioned open mapping results. Section 4 provides several examples, as well as some results for $QDQ$s of curves. Moreover a class of $QDQ$s for set-valued maps taking values on ``$F$-abundant" (see Def. \ref{verydense}) subsets is illustrated. In Section 5 we show that  a {\it set-valued Lie bracket} $[f,g]_{set}(q)$ of Lipschitz continuous  vector fields $f,g$ (see \cite{RampSuss2001},\cite{RampSuss2006},\cite{Ramp2007} and \cite{RampFeleqi}) is a $QDQ$ for the commutator-like multi-flow of $f$ and $g$ at $q$.
\section{DEFINITIONS AND BASIC PROPERTIES}
\subsection{Quasi-Differential Quotients}

	\begin{definition}\label{qdq}	Assume that $F : \rr^n \rightsquigarrow \rr^m$       is a 
		set-valued map, $(\bar x,\bar y) \in \rr^n\times\rr^m  $,   $\Lambda\subset \Lin( \rr^n, \rr^m) $   is a compact set,  and $\Gamma\subseteq\rr^n$  is any  subset.
		We
		say that $\Lambda$ is a {\em Quasi Differential Quotient  ($QDQ$) of $F$ at  $(\bar x,\bar y)$ in the direction of $\Gamma$ }  if there
		exists a $\delta^* \in (0,+\infty]$ and modulus $\rho:[0,+\infty[\to [0,+\infty[$ \footnote{We call{ \it modulus }  a non-decreasing non-negative function $\rho:[0,+\infty[\to [0,+\infty[$ such that $\lim\limits_{\delta\to 0^+}\rho(\delta)=0$.} having the property that
		\begin{itemize}
			\item[$(*)$]
			for every $\delta \in (0,\delta^*)$  there is a continuous map  
			$(L_\delta,h_\delta):\left(\bar{x}+B_{\delta}\right)\cap\Gamma\to \Lin( \rr^n, \rr^m) \times \rr^m$ 
			% with convex compact values 
			such that,	whenever $x
			\in  (\bar x+B_\delta)\cap\Gamma$, 
			\begin{equation} \label{approssimazione}\begin{array}{c}d(L_\delta(x), \Lambda)\leq \rho(\delta), 
			\quad |h_\delta(x)|\leq \delta \rho(\delta)\\\\
				\bar y +  L_\delta(x)\cdot(x-\bar x)  + h_\delta(x)\in F(x) .\end{array}\end{equation}
		\end{itemize}
\end{definition}

The definition of $QDQ$ can be extended to  (set-valued) functions from a manifold  into another  manifold.
\begin{definition}\label{agdqM}  
	Let $\mathcal{N}$, $\mathcal{M}$ be differential  manifolds of class  $C^1$.  Assume that $\tilde G : \mathcal{N} \rightsquigarrow \M$       is a 
	set-valued map, $(\bar x,\bar y) \in \mathcal{N}\times\M  $,   $\tilde\Lambda\subset \Lin( T_{\bar x}\mathcal{N}, T_{\bar y}\M) $   is a compact set,  and $\tilde\Gamma\subseteq\mathcal{N}$  is any  subset.
	Moreover, let $\phi:U\to \rr^n$ and $\psi:V\to\rr^m$ be charts defined on  neighborhoods  $U$ and $V$ of $\bar x$ and $\bar{y}$, respectively, and assume that $\phi(\bar x) = 0$, $\psi(\bar y)=0$. Consider the map 
	$G:=\psi\circ \tilde G\circ\phi^{-1}: \phi(U)\to\rr^m$ and extend it arbitrarily to a map $G:\rr^n\to\rr^m$ (we do not relabel).
	We say that $\tilde\Lambda$ is   a {\em Quasi Differential Quotient ($QDQ$) of $\tilde G$ at  $(\bar x,\bar y)$  in the direction of $\tilde\Gamma$}  if $\Lambda:= D\psi(\bar y) \circ \tilde{\Lambda}\circ D\phi^{-1}(0)$ is a  Quasi Differential Quotient   of $G$ at  $(0,0)$ in the direction of $\Gamma:=\phi(\tilde\Gamma\cap U)$ \footnote{One can easily verify that this definition is indeed chart-independent.} .
\end{definition}

\subsection{Basic properties} We now present some properties which make $QDQ$s fit for proving Maximum Principles via set-separation techniques.
\begin{proposition} \label{loc} Let $\mathcal{M} , \mathcal{N}$ be $C^1$ real differential  manifolds of dimension $m$ and $n$, respectively, and let 
	$F,G:\mathcal{N} \rightsquigarrow \mathcal{M}$
%	$F,G:\rr^n \rightsquigarrow \rr^m$%
	 be set-valued maps. Assume that  $\bar x \in \mathcal{N}$, $\bar y, \bar y_F, \bar y_G \in \mathcal{M}$, $\Gamma, \Gamma_F, \Gamma_G \subseteq  \mathcal{M}$,
	 % $L \in \Lin(\rr^n,\rr^m)$,
	   and $\alpha,\beta \in \rr$. Then: 
	\begin{enumerate}
		\item \textit{[Locality]} If $U$ is a neighborhood of $\bar x$ and $F(x)=G(x)$ for $x \in U\cap \Gamma$, then $\Lambda$ is a $QDQ$ for $F$ at $(\bar x, \bar y)$ in the direction of $\Gamma$ if and only if it is a $QDQ$ for $G$ at $(\bar x, \bar y)$ in the direction of $\Gamma$.
		\item \textit{[Linearity]} If $\mathcal{M}=\rr^m$ and  $\Lambda_F$ and $\Lambda_G$ are $QDQ$ for $F$ and $G$ at points $(\bar x, \bar y_F)$ and $(\bar x, \bar y_G)$ and in the direction of $\Gamma_F$ and $\Gamma_G$, respectively, then $\alpha\Lambda_F+\beta\Lambda_G$ is a $QDQ$ for $\alpha F+\beta G$ at point $(\bar x, \alpha\bar y_F+ \beta\bar y_G)$ in the direction of $\Gamma_F \cap \Gamma_G$.
		\item \textit{[Set product property]} Under the same assumptions as in 2), except for $\mathcal{M}$ which here  is allowed to be any manifold, $\Lambda_F\times \Lambda_G$ is a $QDQ$ at $\Big(\bar x, (\bar y_F, \bar y_G)\Big)$, in the direction of $\Gamma_F\times \Gamma_G$ for the set-valued map $F \times G: x \mapsto F(x)\times G(x)$.
		\item \textit{[Product Rule]} If $\mathcal{M}=\rr$, and still  using the same notation as in $2)$, we have $F(\bar x)\Lambda_G+ G(\bar x)\Lambda_F$ is a $QDQ$ for $FG: x \mapsto F(x)G(x)$.
		\item \label{casospeciale} If $F$ is single-valued   and $L\in  \Lin( T_{\bar x}\mathcal{N}, T_{\bar y}\M) $, $\{L\}$ is a $QDQ$ for $F$ at $(\bar x, \bar y)$ in the direction of $\mathcal{M}$ if and only if $F$ is differentiable at $\bar x$ and $L=DF(\bar x)$.
	\end{enumerate}
\end{proposition}
\noindent \begin{proof} The first property, namely locality, is straightforward from the definition.
	To prove the  linearity stated in $2)$ let   $\rho^F$, $L_\delta^F$,  $h_\delta^F$, and $\rho^G$,  $L_\delta^G$,  $h_\delta^G$ be the maps verifying $(*)$ in Definition \ref{qdq} for  $\Lambda_F$ and $\Lambda_G$, respectively. Hence the thesis follows upon defining  $h_\delta:=\alpha h_\delta^F+\beta h_\delta^G$, $L_\delta:=\alpha L_\delta^F+\beta L_\delta^G$ and $\rho:=|\alpha| \rho^F+|\beta| \rho^G$ .\\
Property $3)$ is straightforward. We  omit the proof  of  $4)$ as well, for it can be trivially obtained by using the same notation as in $2)$, and  multiplying term by term the relations coming from $(*)$ in Definition \ref{qdq}.\\
 Let us prove 5) in the case when  $\mathcal{N}=\rr^m, \mathcal{M}=\rr^n$, the extension of the property to manifolds being trivial. The sufficiency of the differentiability in order $\{L\}$ be a $QDQ$ at $(\bar x, F(\bar x))$ in the direction of $\mathcal{M}$  is trivial.  To prove that  it is also  necessary, observe that from $L_\delta(x)(x-\bar x)=F(x)-F(\bar x)-h_\delta(x)$ and $|L_\delta(x)-L|<\rho(\delta)$, it follows that $|F(x)-F(\bar x)-L\cdot (x-x)|=$ $|(L_\delta(x)-L)\cdot(x-\bar x)+ h_\delta(x)|\le 2\delta\rho(\delta)$ holds for all $\delta>0$, provided $|x-\bar x|\le \delta$. Hence $F$ is differentiable  at $\bar x$ and $L=DF(\bar x)$.
%	
%	{\color{red} at a point $x_0$, one can take $L_\delta(x)\equiv L$ and $h_\delta(x)=f(x)-f(x_0)-\alpha\cdot(x-x_0)$ to have $|h_\delta(x)|\le |x-x_0|\rho(|x-x_0|)<\delta\rho(\delta)$ for a suitable modulus $\rho$ and trivially $f(x_0)+\alpha\cdot(x-x_0)+f(x)-f(x_0)-\alpha\cdot(x-x_0)=f(x)$.}
\end{proof}
\begin{proposition}[Chain rule]\label{chain}
	Let  $\mathcal{N}$ ,$\mathcal{M}$, $\mathcal{L}$ be $C^1$ manifolds of dimensions  $m,n,l$, respectively, let   $F:\mathcal{N} \rightsquigarrow \mathcal{M}$ and $G:\mathcal{M} \rightsquigarrow \mathcal{L}$ be set-valued maps, and consider the composition $G\circ F:x \in \mathcal{N}\rightsquigarrow \bigcup\limits_{y \in F(x)} G(y) \in\mathcal{L}$.\\
	Assume $\Lambda_F$ is a $QDQ$ for $F$ at point $(\bar x, \bar y)$ in the direction of $\Gamma_F$ and $\Lambda_G$ is a $QDQ$ for $G$ at point $(\bar y, \bar z)$ in any direction $\Gamma_G$ containing $F(\Gamma_F)$. Then the set $\Lambda:=\Lambda_G \circ \Lambda_F$ of all compositions of elements of $\Lambda_G$ with elements of $\Lambda_F$ is a $QDQ$ for $G \circ F$ at $(\bar x, \bar z)$ in the direction of $\Gamma_F$.
	
\end{proposition}
\begin{proof}
Let us prove Proposition \ref{chain} when 	$\mathcal{M}=\rr^m$, $\mathcal{N}=\rr^n$, and $\mathcal{L}=\rr^l$, the locality and chart independence  of the notion of $QDQ$ guaranteeing that the proof is still valid in the general case
	
	Let $\delta_*^F$, $\delta_*^G$,   $L_\delta^F$, $L_\delta^G$, $h_\delta^F$, $h_\delta^G$, $\rho^F$, $\rho^G$, with $0<\delta<\delta_*:=\min\{\delta_*^F,\delta_*^G\}$,  the maps  involved in the definition of the $QDQ$s $\Gamma_F$ and $\Lambda_G$.  In particular, one has, for all $\delta\in(0,\delta_*)$, 
	\begin{equation}\label{one}\bar y +  L^F_\delta(x)\cdot(x-\bar x)  + h^F_\delta(x)\in F(x), \, \forall x \in (\bar x+B_\delta)\cap \Gamma_F, \end{equation}
	\begin{equation}\label{two}\bar z + L^G_\delta(y)\cdot (y-\bar y)+ h_\delta^G(y) \in G(y), \, \forall y \in (\bar y+B_\delta)\cap F(\Gamma_F).\end{equation}
	Let us set $$M:=\max\left\{\max\limits_{L' \in \Lambda_F}|L'|,\max\limits_{L'' \in \Lambda_G}|L''|,1\right\}, \quad \eta(\delta):=\frac{\delta}{3M},\,\,\forall \delta\geq 0.$$ If needed, let us redefine  $\delta^*$  in order that   $\delta^*<1$ and it  small enough to guarantee $\rho^F\left(\eta(\delta)\right)\le M$ for all $\delta<\delta^*$.\\ 
	Now, let us consider any $x \in (\bar x+B_{\eta(\delta)})\cap \Gamma_F$ and let us set  $y:=\bar{y}+L_{\eta(\delta)}^F(x)\cdot(x-\bar x)+h^F_{\eta(\delta)}(x)$. Clearly, we have $y \in F(\Gamma_F)$, and  $|y-\bar y|\leq {\eta(\delta)}[M+2\rho^F({\eta(\delta)})]\le \delta$, so that $y \in (\bar y+ B_\delta)\cap F(\Gamma_F)$.\\
	If we set $\xi(x):=\bar{y}+L_{\eta(\delta)}^F(x)\cdot(x-\bar x)+h^F_{\eta(\delta)}(x)$,  we get  
%	\begin{multline}\label{three}\bar z + L^G_\delta\Bigg[\bar{y}+L_{{\eta(\delta)}(\delta)}^F(x)\cdot(x-\bar x)+h^F_{\eta(\delta)}(x)\Bigg]\cdot \Big(L_{\eta(\delta)}^F(x)\cdot(x-\bar x)+h_{\eta(\delta)}^F(x)\Big)+\\+ h_\delta^G\Big[\bar{y}+L_{\eta(\delta)}^F(x)\cdot(x-\bar x)+h^F_{\eta(\delta)}(x)\Big] \in G(F(x)),\end{multline}
%	holds $\forall x \in (\bar x+B_{\eta(\delta)})\cap \Gamma_F$. This can be rewritten as
\begin{multline}\label{four}
	\bar z+ \Bigg\{L^G_\delta(\xi(x))\circ L_{\eta(\delta)}^F(x)\Bigg\} \cdot(x-\bar x)+{h}_{\eta(\delta)}(x) \in G(F(x)) \\
\forall x \in (\bar x+B_{\eta(\delta)})\cap \Gamma_F ,
	\end{multline}
	where the continuous function ${h}_{\eta(\delta)}$, defined as $${h}_{\eta(\delta)}(x):=L^G_\delta(\xi(x))\cdot h_{\eta(\delta)}^F(x)+ h_\delta^G(\xi(x))$$  verifies \bel{stimah}|{h}_{\eta(\delta)}(x)|\le M{\eta(\delta)}\rho^F({\eta(\delta)})+3M{\eta(\delta)}\rho^G(3M{\eta(\delta)}).\eeq
The function between curly brackets in  \eqref{four}, i.e. $$\delta\mapsto L_{\eta(\delta)}(x):=L^G_\delta\Big(\xi(x)\Big)\circ L_{\eta(\delta)}^F(x),$$  is also a continuous function as it is composition of continuous functions.
	Moreover, every image  $L_{\eta(\delta)}(x)$ is the composition of a linear functional whose distance from $\Lambda_G$ is at most $\rho^G(3M{\eta(\delta)})$ and whose norm is at most $M$ with a linear functional whose distance from $\Lambda_F$ is at most $\rho^F({\eta(\delta)})$. Therefore \bel{stimagamma} d\left(L_{\eta(\delta)}(x),  \Lambda\right) \leq M[\rho^F({\eta(\delta)})+\rho^G(3M{\eta(\delta)})].\eeq
Finally, if one  sets  
 $$ \rho({\eta(\delta)}) :=M\Big[\rho^F({\eta(\delta)})+3\rho^G(3M{\eta(\delta)})\Big]$$
 then formulas \eqref{four}-\eqref{stimagamma} tell us that $\Lambda$ is a $QDQ$ for $G \circ F$ at  $(\bar x,\bar z)$ in the direction of $\Gamma_F$. 
	
\end{proof}

 \section{Open mapping and set-separation} 
  
 \subsection{Open mapping}
  Open mapping results are essential for any reasonable  generalized differentiation theory.  We recall here a  result  which, in particular, marks the difference between $QDQ$s and $AGDQ$s.
   \begin{theorem}[Open Mapping for QDQs]
 	Let $F:\rr^n \rightsquigarrow \rr^m$ be a set-valued map and let $\Gamma$ be a convex cone in $\rr^n$. Let $\Lambda$ be a $QDQ$ of $F$ at $(\bar x, \bar y)$ in the direction of $\Gamma$ and assume that $L\cdot \Gamma=\rr^m$ for all $L$ in $\Lambda$. Then:
 	\begin{enumerate}
 		\item[$(i)$] There exist constants $\alpha,\beta>0$ such that $\forall a \in (0,\alpha]$ $$\bar y + (B_a \setminus \{0\}) \subseteq F(\bar x +B_{a\beta}\cap \Gamma).$$
 		\item[$(ii)$] There exist $\delta^*>0$ such that for all $\delta<\delta^*$ and every $(L_\delta,h_\delta)$ as in Definition \ref{qdq}, there exists $x_\delta \in (\bar x+ B_\delta) \cap \Gamma$ such that:
 		$$\bar y=L_\delta(x_\delta)(x_\delta-\bar x)+h_\delta(x_\delta)  \,\, \left(\in F(x_\delta)\right).$$
 	\end{enumerate}
 	In particular, by possibly reducing $\alpha$, the following inclusions hold for all $a \in (0,\alpha]$
 	\bel{nonp}\bar y + B_a \subseteq F(\bar x +B_{a\beta}\cap \Gamma).\eeq
 \end{theorem}
 \begin{remark}
	Part (ii)  is not true if $\Lambda$ is just a $AGDQ$. In that case, one can only prove  a ``punctured'' inclusion, namely in general $ \bar y\notin  F(\bar x +B_{a\beta}\cap \Gamma).$	
	\end{remark}
 \subsection{Transversal cones}
 
 Let us recall the notions of transversality and strong transversality of cones. Let $E$ be a finite-dimensional, real  linear space, and let $E^*$ be  its dual space.  A subset $\K\subseteq E$   is a {\it cone} if $\alpha k\in \K$  for all $(\alpha,k)\in [0,+\infty[\times \K$. If $D\subseteq E$ is any subset, let us set 
 $$
 \begin{array}{l}
 	\Span^+D :=  \ds\Big\{\sum_{i=1}^\ell \alpha_i v_i: \  \ell\in \mathbb{N},\ \alpha_i\geq 0, \,\,v_i\in D,\,\,  \forall i\le l  \Big\}, \\ 
 	D^{\bot} :=\Big\{p\in E^* :\   p\cdot w \leq 0 \ \ \forall \;
 	w\in D \Big\}\subseteq E^*.
 \end{array}
 $$
 The convex cones $\Span^+D$ and $D^{\bot}$ are called the {\em conic hull} of $D$ and 
 the {\em polar cone} of $D$, respectively.  Let   $\K_1$, $\K_2\subseteq E $ be convex cones. We say that 
 $\K_1$ and $\K_2\subseteq E$ are {\em transversal}, if $
 \K_1-\K_2:=\big\{k_1-k_2 :\ (k_1,k_2)\in  \K_1\times \K_2\big\} = E$. 
 $\K_1$ and $\K_2$  are {\em strongly transversal},
 if they are transversal  and $\K_{1}\cap\K_{2} \supsetneq\{0\}$.
 
 \begin{proposition}\label{teo2} Two  convex cones $\K_1$, $\K_2\subseteq E $ are transversal if and only if they are either strongly transversal or complementary linear subspaces, namely  $\K_1\oplus{\K}_2=E$  (i.e., $\K_1 +\K_2=E$ and $\K_1 \cap\K_2=\{0\}$).
 \end{proposition}
 
 Saying that two cones {\rm$\K_1$ and $\K_2$ are not strongly transversal} is equivalent to saying that they are linearly separable:
 \begin{proposition}\label{teo3}Two  convex cones $\K_1,\K_2\subseteq E$ are not transversal  if and only if \rm$\K_1$ and $\K_2$ are  {\rm linearly separable}, by which we mean that   $(-\K_1^\bot\cap\K_2^\bot) \backslash \{0\} \neq \emptyset,$ namely, there exists a  linear form $\lambda\in E^*\backslash \{0\}$
 	such that, $\forall (k_1,k_2)\in \K_1\times  \K_2$, 
 	$ \lambda \cdot k_1\geq 0$, 
 	and $\lambda \cdot k_2\leq 0$.
 \end{proposition} 
 In  a linear space  $E$, let us call {\it convex multi-cone} any family of convex cones of $E$. 
 
 \begin{definition}{\rm \cite[Def. 2.5]{PR}}\label{ApprCone}  Let $\M$ be a $C^1$ differentiable manifold, $\mathcal E\subseteq\N$   a set, and $z\in \mathcal E$.  A  {\em  $QDQ$ approximating multi-cone}  to $\mathcal E$ at $z$ is a
 	convex multi-cone $\mathbf{K}\subseteq T_{z}\N$ such that there exist an integer   $n\ge0$, a set-valued
 	map  $G : \rr^n  \rightsquigarrow \N$, a convex cone $\Gamma\subseteq\rr^n$, and a Quasi Differential Quotient  $\Lambda$ of $G$  at  $(0,z)$ in the direction of $\Gamma$ such that $G(\Gamma)\subseteq \mathcal E$ and $\mathbf{K} =
 	\{L\cdot\Gamma\  : L\in\Lambda\}$.\footnote{When    a
 		$QDQ$ approximating multi-cone  is a singleton, namely $\mathbf{K}=\{K\}$, we say  that $K$ is    {\em  a
 			$QDQ$ approximating cone}  to $\mathcal E$ at $z$.}
 	We say that such a {\em triple $(G,\Gamma,\Lambda)$  generates the multi-cone $\mathbf{K}$}.
 	
 	If the triple $(G,\Gamma,\Lambda)$ defining a $QDQ$ approximating multi-cone $\mathbf{K}$   can be chosen so that 
 	$G(\Gamma)\subseteq \mathcal E\backslash \{z\}$,  then we say that the   $\mathbf{K}$
 	is {\rm $z$-ignoring}. 
 \end{definition}
 
 \begin{remark}\label{Rignoring}{\rm Because of  local character of the notion of $QDQ$ for a set-valued map, one can equivalently say that {\it a $QDQ$ approximating multi-cone  $\mathbf{K}$ is $z$-ignoring if  there is some $\delta>0$ such that $G(B_\delta\cap\Gamma )\subseteq \mathcal E\backslash \{z\}$ for some $\delta>0$.}}\end{remark}
 \begin{remark}{\rm For single-valued maps, $QDQ$ approximating multi-cones consisting of a single cone are precisely Boltyanski's approximating cones.}\end{remark}
 
 \begin{definition} Let $\mathcal{X}$ be a topological space, and let   $\mathcal{D}_1, \mathcal{D}_2\subseteq\mathcal{X} $, $y\in  \mathcal{D}_1\cap \mathcal{D}_2$. We say that $\mathcal{D}_1$ and $ \mathcal{D}_2$ are {\em  locally separated at $y$} provided there exists a neighborhood $V$ of $y$ such that
 	$
 	\mathcal{D}_1\cap\mathcal{D}_2\cap V = \{y\}.
 	$
 \end{definition}
 The following open-mapping-based  result, obtained in \cite{PR},  characterizes  set separation in terms of linear separation of  $QDQ$ approximating cones.  It includes a crucial approximation property (see  (ii)  below)  whenever one of the cones is $z$-ignoring.
 \begin{theorem}{\rm \cite[Thm. 2.3]{PR}}\label{teoteo}
 	Let  $ {\mathcal E}_1, {\mathcal E}_2 $ be subsets of a $C^1$ differentiable manifold  $\N $  and 
 	let
 	$ z\in {\mathcal E}_1\cap {\mathcal E}_2 $.
 	Assume that  $\K_1$, $\K_2$ are $QDQ$ approximating cones  of  $ {\mathcal E}_1 $ and $ {\mathcal E}_2 $, respectively, at $z$. 
 	\begin{itemize}
 		\item[\em (i)] If  $\K_1$ and $\K_2$ are strongly transversal,
 		then  the sets ${\mathcal E}_1$ and ${\mathcal E}_2$
 		are not locally separated. 
 		
 		\item[\em (ii)] If $\K_1$ (or, equivalently, $\K_2$) is $z$-ignoring and
 		$\K_1$ and $\K_2$ are transversal then the sets ${\mathcal E}_1$ and ${\mathcal E}_2$
 		are not locally separated.
 	\end{itemize}
 \end{theorem}

\section{Examples}
\begin{example}[The $\delta$-independent case]\label{exnodelta}Let us begin observing that if, for a given set-valued map $F : \rr^n \rightsquigarrow \rr^m$,   a point $(\bar x,\bar y) \in \rr^n\times\rr^m  $, a subset  $\Gamma\subseteq\rr^n$  $\Lambda\subset \Lin( \rr^n, \rr^m) $   is a compact set such that 
	$$ \displaystyle \lim_{\delta\to 0}\|d(L(x), \Lambda)\|_{L^\infty(\bar x + B_\delta \cap \Gamma)}=0, 
	\quad {\|h(x)\|_{L^\infty(\bar x + B_\delta \cap \Gamma)}} = o(\delta)$$and
	$	\bar y +  L(x)\cdot(x-\bar x)  + h(x)\in F(x),$
	where $(L,h): \Gamma \to \Lin(\rr^n,\rr^m)\times \rr^m$ is a ($\delta$-independent) continuous function, then   $\Lambda$ is a $QDQ$. Indeed, one can check that,  if  $(L_\delta,h_\delta):=(L,h)_{\bar x + B_\delta \cap \Gamma}$  and  $$\rho(\delta):=\max\left\{\|d(L(x), \Lambda)\|_{L^\infty(\bar x + B_\delta \cap \Gamma)},\frac{\|h(x)\|_{L^\infty(\bar x + B_\delta \cap \Gamma)}}{\delta} \right\},$$ for   every $\delta>0$, then the properties in Definition \ref{qdq} are met. 
	
	However it is not, in general, possible to reduce $\delta\mapsto(L_\delta,h_\delta)$ to a constant in the definition of $QDQ$, as one can check already in the $n=m=1$ case below.
	
\end{example}
\begin{example}[Clarke's generalized  Jacobian]  Let $\Omega\subseteq\rr^n$ be  an open set, let $F:\Omega\to\rr^m$  be a (single-valued) Lipschitz continuous map, and let    $x^*\in\Omega$. In \cite{nostropaper}, we proved a sufficient condition for a set to be a $QDQ$ via approximation of the function: it is possible to infer from that the { Clarke's generalized  Jacobian}\footnote{ Let us recall that  the  {\it Clarke's generalized  Jacobian} $\partial_CF(\bar x )$ at $\bar x$ is defined as the set  $\overline{\text{co}} \left\{v = \displaystyle\lim_{n\to \infty} DF(x_n), \quad (x_n)_{n\in\nn}\subset \text{Diff}(F), \,x_n\to x^*\right\},$ 
		where $\text{Diff}(f)$  denotes the set comprising all points $x\in\Omega$ such that $F$ is differentiable at $x$.  
	} $\partial_CF(\bar x)$ at $\bar x$  is a $QDQ$ for $F$ at  $(\bar x,F(\bar x))$ in the direction of $\rr^n$.
\end{example}

\subsection{The case $n=1$}	

\begin{example} \label{exampleabsolutevalue1d} Since $[-1,1]$ is the Clarke's generalized Jacobian of the function $f(x):=|x|$ at $x=0$, in view of the previous example it is a $QDQ$ for the function $f$ at $(0,0)$ in the direction of $\rr$.
	However, we prefer here to give an explicit construction of the maps $L_\delta, h_\delta$. 
	For every $\delta>0$, let us set 
		$$
		g_{\delta}(x) := \left\{ \begin{array}{ll}\displaystyle\frac 12 \left(\frac{x}{{\delta}}\right)^2 + \frac{{\delta}^2}{2}\quad&\forall x\in[-{\delta}^2, {\delta}^2]\\\\
			|x|&\forall   x\in\rr\backslash [-{\delta}^2, {\delta}^2].\end{array}\right.
		$$
		%	One has, 
		%	$$g_{\delta}(-{\delta}^2) = g_{\delta}({\delta}^2) = {\delta}^2 = f({\delta}^2) = f(-{\delta}^2)
		%	$$so that $g_{\delta}$ is continuous. Moreover,
		%	$$
		%	g_{\delta}'(x) = \left\{ \begin{array}{ll}\displaystyle \frac{x}{{\delta}^2}\quad&\forall x\in[-{\delta}^2, {\delta}^2]\\\\
		%	\sgn(x)&\forall   x\in\rr\backslash [-{\delta}^2, {\delta}^2]\end{array}\right.
		%	$$
		%	hence, 
		For every ${\delta}>0$, $g_{\delta}$ is of class $C^1$.
		%19748
		Moreover, for  every $x\in[-\delta,\delta]$, one has 
		$
		f(x)=|x|=
		L_{\delta}(x)\cdot  x + h_{\delta}(x), 
		$
		where 
		$
		L_{\delta}(x):= \frac{x}{{\delta}^2}{\bf 1}_{[-{\delta}^2,{\delta}^2]} +  \sgn(x){\bf 1}_{[-{\delta},-{\delta}^2]\cup [{\delta}^2,{\delta}]}
		$
		and
		$
		h_{\delta}(x) := \frac{{\delta}^2}{2} +   (|x|-g_{\delta}(x)){\bf 1}_{[-{\delta}^2,{\delta}^2]}.
		$
		Notice that  $L_{\delta}([-{\delta},{\delta}]) = [-1,1]$, and $|h_{\delta}(x)|\leq {\delta}^2$, so $\Lambda= [-1,1]$ meets Definition \ref{qdq} with $\rho(\delta)=\delta$.

		\vskip0.7truecm
	
		By definition of  $QDQ$ any compact set $\Lambda'\supseteq \Lambda$ is again a $QDQ$ for the same $f$ at the same point $(0,0)$ in the same direction. Hence, in this example  $[-1,1]$ is  the smallest possible $QDQ$ for $f$ at $(0,0)$ in the  direction of $\mathbb{R}$.
				\end{example}
			 Actually, while the previous example is a particular case of Clarke's generalized Jacobian (see below), it is also is an instance of the following general fact:
%		
%		{\fra  NELLA PROP. SEGUENTE  MI PARE CHE sia inutile fare nella direzione di un Neighborhood, a quel punto è lo stesso ed è più chiaro dire nella direzione di $\rr$ }
		\begin{proposition} \label{observation}
			Let $\mathcal{M}$ be an $m$-dimensional manifold of class $C^1$ and $f:\rr \to \mathcal{M}$
			 be a continuous curve admitting left and right derivatives  $f'(\bar t^-)$ and $f'(\bar t^+)$ at a point  $\bar t \in\rr$. Any compact set $\Lambda \subseteq \Lin(\rr,T_{f(\bar t)}M)$ containing an arc connecting $f'(\bar t^+)$ and $f'(\bar t^-)$ is a $QDQ$ for $f$ at $(\bar t, f(\bar t))$ in  the direction of  $\rr$. Conversely, any $QDQ$ for $f$ at $(\bar t, f(\bar t))$ in the direction of $\rr$ necessarily contains a connected set containing $\{f'(\bar t^+),f'(\bar t^-)\}$.
		\end{proposition}
		\begin{proof}
			We provide the proof of the result for $\mathcal{M}=\rr^m$, the proof for the general case being a straightforward consequence which can be obtained by means of a chart. Assume $\bar t=0$ and $f(\bar t)=\mathbf{0}$, which is without loss of generality by linearity and chain rule, and assume there is a curve $\gamma(t)$ defined on $[-1,1]$ and connecting $f'(0^+)$ and $f'(0^-)$ whose support is contained in $\Lambda$. Let us also name $s_{\pm}(t)$ the segments defined on $\left[-\delta^2,\frac{-\delta^2}{2}\right]$ and $\left[\frac{\delta^2}{2},\delta^2\right]$ respectively, connecting $\frac{f(-\delta^2)}{- \delta^2}$ with $f'(0^{-})$ and $f'(0^{+})$ with $\frac{f(-\delta^2)}{- \delta^2}$, respectively. We define, for every $\delta \in (0,1)$: $$L_\delta(t):=\begin{cases}
			\frac{f(t)}{t} & \text{ if } |t| \in [\delta^2,\delta)\\
			\gamma\left(\frac{2t}{\delta^2}\right) & \text{ if } |t| \le \frac{\delta^2}{2} \\
			s_{\pm}(t) & \text{ if } |t| \in \left[\frac{\delta^2}{2},\delta^2\right].
			\end{cases}$$
			It is clear that $L_\delta(t)$ was constructed continuous with a proper choice of the affine functions $s_{\pm}$ on the two intervals $\left[-\delta^2,-\frac{\delta^2}{2}\right]$ and $\left[\frac{\delta^2}{2},\delta^2\right]$. Now, $L_\delta\left(\left[-\frac{\delta^2}{2},\frac{\delta^2}{2}\right]\right)$ was constructed to be contained in $\Lambda$ and $L_\delta\left([-\delta,\delta]\setminus \left[-\frac{\delta^2}{2},\frac{\delta^2}{2}\right]\right)$ consists of points whose distance from $f'(\bar t^+)$ and $f'(\bar t^-)$ is bounded by a suitable modulus $\hat{\rho}(\delta)$ by definition of left and right derivative as limit of one-sided differential quotients.\\
			At the same time, if we define $h_\delta(t)$ to be $$h_\delta(t):=f(t)-L_\delta(t)t$$ then it is by definition continuous and $h_\delta(t)$ is clearly $0$ outside $[-\delta^2,\delta^2]$ and at $t=0$. For $0<|t|<\delta^2$, let us consider $|h_\delta(t)|/\delta^2 \le \frac{|f(t)|}{\delta^2}+|L_\delta(t)| \le M$ for a suitable constant $M$. Choosing $\rho(\delta)=\max\{M\delta,\hat{\rho}(\delta)\}$ we have all the properties $(*)$ satisfied hence $\Lambda$, which was an arbitrary compact set containing a curve connecting the left and right derivative, is a $QDQ$ for our function at $(0,0)$, in the direction of $\rr$.\\ To prove the converse, let $\Lambda$ be a $QDQ$ for $f$ at $(0,0)$ in the direction of $\rr$ and notice that, by the definition of $QDQ$,  the points $$L_\delta(\pm\delta)=\frac{f(\pm \delta)}{\pm \delta}-\frac{h_\delta(\pm \delta)}{\pm\delta}$$ have  vanishing distances from the compact set $\Lambda$ as $\delta \to 0$.  At the same time, they  converge  to $f'(0^{\pm})$, respectively. Hence $\{ f'(0^-), f'(0^+)\}\subseteq\Lambda$. On the other hand let us assume by contradiction that $\Lambda$ can be split into two disjoint compact sets $C_-$ and $C_+$ containing $f'(0^-)$ and $f'(0^+)$ respectively. The distance function $d(x,y): C_-\times C_+ \to \rr$ is a continuous function defined on a compact set, hence it has a minimum $\varepsilon$ and $\varepsilon >0$ and the sets $C_{\pm}^\varepsilon=\left\{L \in \Lin(\rr,\rr^m): d(L,C_{\pm})\le \frac{\varepsilon}{3}\right\}$ are still disjoint. If $\delta$ is small enough to guarantee $\rho(\delta)\le \frac{\varepsilon}{3}$ and $|L_\delta(\pm\delta)-f'(0^{\pm})|\le \frac{\varepsilon}{3}$, the connected set $L_\delta([-\delta,\delta])$, consisting only of points whose distance from $\Lambda$ is smaller than $\rho(\delta)$, is necessarily wholly contained in just one of the $C_{\pm}^\varepsilon$, which is in contradiction with the fact that $L_\delta(\pm\delta) \in C_{\pm}^\varepsilon$.
		\end{proof}
	\begin{remark}
		As a  consequence of the above Proposition,  if $\mathcal{M}$ is a $1$-dimensional $C^1$ manifold  and if  $f:\rr \to \mathcal{M}$
		is a continuous curve admitting left and right derivatives  $f'(\bar t^-)$ and $f'(\bar t^+)$ at a point  $\bar t \in\rr$, then {\it the segment $$\Big[\min\{f'(\bar t^-),f'(\bar t^+)\}, \max\{f'(\bar t^-),f'(\bar t^+)\} \Big]$$ is the smallest $QDQ$ of $f$ at $t$ in the direction of $\rr$.} However, as soon as  $m>1$, the  $QDQ$s of curves, even those which are minimal with respect to inclusion, may very well be non-convex subsets (but always containing a connected set).
	\end{remark}
%\textcolor{purple}{GNU (legga e poi cancelliamo): La ragione per cui questo argomento non può essere replicato in favore della connessione per archi è che le componenti connesse per archi di un compatto non sono necessariamente oggetti a distanza positiva, come un famoso esempio con $\{\left(x,\sin\left(\frac 1x\right)\right), \, x \in (0,1)\}\cup \{0\}\times [-1,1]$ dimostra. Il meglio che si può fare con questa tecnica è "Se un $QDQ$ è localmente connesso per archi, contiene un arco che collega le derivate sinistre e destre." ------------------------------------------ Il lettore capisce immediatamente che questa proposizione diventa una caratterizzazione a tutti gli effetti nel caso $n=m=1$?}
%	
%

\subsection{Set-valued functions ranging in ``$F$-abundant" subsets}
Let us begin with the definition of {\it $F$-abundant subset}:
\begin{definition} \label{verydense} Let $\E\subseteq E\subseteq \rr^m$ be subsets  and let  $F:\rr^n\to E\subseteq \rr^m$ be a (single-valued) map.  We say that $\E$ is  $F$-{\it abundant} if,  for every $\eta>0$, there exists a continuous  map $\theta_\eta:\rr^m\to\E$ such that
	$|F(x) - \theta_\eta\circ F(x)|<\eta, $ for all $x\in \rr^n$. \footnote{The notion of  {$F$-abundant set } generalizes the control-theoretical concept of {\it abundant set}, which was introduced by J. Warga \cite{warga1} and successively elaborated by B. Kaskosz \cite{Kaskosz}.}
\end{definition}
We illustrate here    a situation where  a known $QDQ$  $\Lambda$ of  a map  $F:\rr^n\to E\subseteq \rr^m$ happens to be also a  $QDQ$  for a set-valued map $\tilde F:\rr^n\rightsquigarrow \E\subseteq E$, provided  $\E$ is $F$-{\it abundant}  in $E$.
Such a case is found in applications,  in particular for the {\it infimum gap} problem, e.g. when  approximating the reachable set of a control system with  ($QDQ$) approximating cones to  the larger reachable set of the convexified system (see \cite{PR}). 

\begin{example}
	{\it Let $F : \rr^n \to E\subseteq \rr^m$       be a 
		(single-valued) map, $\bar x \in \rr^n$,  and let  $\Lambda$  be a $QDQ$ of  $F$ at $\left(\bar x,F(\bar x)\right)$ in the direction   of  a given set  $\Gamma\subseteq\rr^n$. 
		If    $\E$ is $F$-{\it abundant} in $E$
		then $\Lambda$  be a $QDQ$ of  the set-valued map $$
		\tilde F : \rr^n \rightsquigarrow \rr^m, \qquad  \tilde F(x): = \bigcup_{\eta>0}  \theta_\eta\circ F(x),
		$$ (where the maps  $\theta_\eta$ are  as in Definition \ref{verydense})  at $\left(\bar x,F(\bar x)\right)$ in the direction   of  $\Gamma\subseteq\rr^n$).
	}
	
	Indeed,  for every $\delta>0$, 
	let $\rho(\delta)$, $L_\delta$, and $h_\delta$ 
	% with convex compact values 
	be as in Definition \ref{qdq}. Defining, for every $x\in\rr^n$  and any $\delta>0$,
	$$
	\tilde L_\delta(x) :=  L_\delta(x), \quad  \tilde h_\delta(x):= \left(\theta_{\delta\rho(\delta)}\circ F(x) - F(x)\right) +  h_\delta(x) 
	$$
	and $\tilde\rho(\delta) := 2\rho(\delta)$, 
	we obtain 
	$$\min_{L'\in\Lambda}|\tilde L_\delta(x) - L'|\leq \tilde\rho(\delta), 
	\quad  |\tilde h_\delta(x)|\leq \delta \tilde\rho(\delta),$$ $$
	\bar y +  \tilde L_\delta(x)\cdot(x-\bar x)  + \tilde h_\delta(x) = \theta_\delta\circ F(x) \in \tilde F(x)$$
	whenever $x
	\in  (\bar x+B_\delta)\cap\Gamma$. Hence $\Lambda$ is  $QDQ$ of $
	\tilde F$
	at $\left(\bar x,F(\bar x)\right)$ in the direction   of  $\Gamma\subseteq\rr^n$.

\end{example}

\section{Lie Brackets of Lipschitz vector fields}
Let $f,g:\Omega \subseteq \rr^n \to \rr^n$ be two Lipschitz vector fields and $q\in \Omega$. For $\eta>0$ sufficiently small, and for $t\in (-\eta,\eta)$, we shall use $\Phi_t^f(q)$ to denote  the value at $t$ of the   solution of $$\begin{cases}
	y'(\tau)=f(y(\tau))\\
	y(0)=q
\end{cases}.$$  As soon as  $f,g \in C^1$, it is well known that the multi-flow  $$\Psi_t(q):=[\Phi_{t}^{-g}\circ\Phi_{t}^{-f}\circ \Phi_t^g\circ \Phi_t^f](q)$$ verifies  $\Psi_t(q) = q +  t^2 [f,g](q) +o(t^2)$, where $[f,g]$ is the Lie bracket, namely  $[f,g]=Dg\cdot f- Df \cdot g$. In \cite{RampSuss2001}, a set-valued generalization $[f,g]_{set}$ of the Lie bracket has been introduced for Lipschitz vector fields. It  is defined  by setting, for every $q$  in  the  common domain of $f$ and $g$,  $$[f,g]_{set}(q)=\overline{\text{co}}\left\{\lim\limits_{n\to \infty} [f,g](q_n)\,\right\},$$ where limits are taken over sequences $q_n\to q, \, q_n \in \text{Diff}(f)\cap \text{Diff}(g)$,  $\text{Diff}(f)$  and $\text{Diff}(g)$ denoting the sets of differentiability points of $f$ and $g$, respectively.\footnote{For any subset $S \in \rr^n$, $\overline{\text{co}}(S)$ denotes the  closure of the convex hull of $S$.} Using this bracket, the authors of  \cite{RampSuss2001}, proved a Chow's  type theorem to non-smooth fields. Furthermore, in \cite{RampSuss2006} it was shown that these multi-valued Lie brackets can be used to test commutativity of non-smooth fields and to estimate $\Psi_t(q)-q$ to the second order, while in \cite{Ramp2007} they were used to extend Frobenius' Theorem. Lastly, in an oncoming article \cite{nostropaper} we  are going to establish  Goh-like second order conditions for optima, by means of these set-valued Lie brackets. For this aim we need the following  result:
\begin{proposition} \label{QDQLie}
	Let $\Omega$ be a compact subset of $\rr^n$, $f,g:\Omega \to \rr^n$ be Lipschitz vector fields and $q \in \mathring{\Omega}$. The set $[f,g]_{set}(q)$ is a $QDQ$ for the function $$F:\varepsilon \mapsto \Psi_{\sqrt{\varepsilon}}(q)$$ at point $(0,q)$ in the  direction of  $\rr^+$.
\end{proposition}
\begin{proof}
We are going to   make use of the following integral  formula, established in \cite{RampSuss2001}, \footnote{More precisely, the authors of \cite{RampSuss2001} proved an {\it exact} integral formula, of which \eqref{exactflow} is a straightforward consequence.}
	\begin{multline} \label{exactflow}
		\qquad\qquad\qquad\Phi_{\sqrt{\varepsilon}}^{-X_2}\circ\Phi_{\sqrt{\varepsilon}}^{-X_1}\circ\Phi_{\sqrt{\varepsilon}}^{X_2}\circ\Phi_{\sqrt{\varepsilon}}^{X_1}(q)=\\=q+\int_0^{\sqrt{\varepsilon}}\int_0^{\sqrt{\varepsilon}} [X_1,X_2](\theta(q,\sigma,{\sqrt{\varepsilon}},\tau))\,d\tau\,d\sigma+ o(\varepsilon),
	\end{multline}
	where $X_1,X_2$ are $C^1$ vector fields, and $\theta(x,\sigma,t,\tau)=\Phi_{\tau-t}^{X_1} \circ \Phi_{\tau}^{X_2}\circ \Phi_t^{X_1}(x).$
	Let us consider a mollifier, i.e. a $C^\infty$ function $\varphi(x)\ge 0 $ having  support contained in $B_1(0)$   and $L^1$ norm equal to $1$. For every $\eta>0$, we set $$F_\eta(\varepsilon):=\Phi_{\sqrt{\varepsilon}}^{-g_\eta}\circ\Phi_{\sqrt{\varepsilon}}^{-f_\eta}\circ \Phi_{\sqrt{\varepsilon}}^{g_\eta}\circ \Phi_{\sqrt{\varepsilon}}^{f_\eta}(q),$$ where  $f_\eta$ and $g_\eta$ denote the convolutions  $f*\varphi_\eta$ and $g*\varphi_\eta$, respectively, with  $\varphi_\eta:=\frac{1}{\eta^n}\varphi\left(\frac{x}{\eta}\right)$. If we take $X_1:=f_\eta, X_2:=g_\eta$ in \eqref{exactflow}, and moreover,   we set $$\eta=\eta(\varepsilon):=\varepsilon^2,$$ we  get
	\begin{multline*}\qquad\qquad\qquad\qquad F_\eta(\varepsilon)-F_\eta(0)=\\=\int_0^{\sqrt{\varepsilon}}\int_0^{\sqrt{\varepsilon}}[f_{\eta(\varepsilon)},g_{\eta(\varepsilon)}](\theta(q,\sigma,\sqrt{\varepsilon},\tau))\,d\tau\,d\sigma+o(\varepsilon).\end{multline*}
	From this inequality, and observing that there exists a constant $c>0$ (depending on the Lipschitz constant and $L^\infty$ norm of $f$ and $g$) such that $$|F_{\eta}(\varepsilon)-F(\varepsilon)|<c{\eta(\varepsilon)} \quad \forall {\eta},\epsilon>0,$$ we get $$
	F(\varepsilon)-F(0)=\varepsilon L(\varepsilon)+o(\varepsilon)$$
%	\footnote{Sometimes in what follows it will be notationally convenient to write ${\eta(\varepsilon)}$, but the reader should remember ${\eta(\varepsilon)}={\eta(\varepsilon)}(\varepsilon)=\varepsilon^2$ eliminates in fact freedom for this mollification parameter.}
  where $$L(\varepsilon):=\frac 1\varepsilon \int_0^{\sqrt{\varepsilon}}\int_0^{\sqrt{\varepsilon}}[f_{\eta(\varepsilon)},g_{\eta(\varepsilon)}](\theta(q,\sigma,\sqrt{\varepsilon},\tau))\,d\tau\,d\sigma.$$ It is clear that $L(\varepsilon)$ depends continuously on $\varepsilon$ and that the same  holds true for the map  $h(\varepsilon):=F(\varepsilon)-F(0)-\varepsilon L(\varepsilon)=o(\varepsilon)$.\\ In view of Example \ref{exnodelta}, to prove that $[f,g]_{set}(q)$ is indeed a $QDQ$ of $F$ at $(0,q)$ in the direction of $\rr^+$, it only remains to check that $L(\varepsilon)$ has vanishing distance from it as $\varepsilon$ approaches $0$.  To this aim, let us set $D:=\text{Diff}(f)\cap \text{Diff}(g)$, and, for every $\varepsilon>0$,  $D^\varepsilon(\theta):=B_1(0)\cap\{v:\, \theta+{\eta(\varepsilon)} v \in D\},$   and
	$$Dg_{\eta(\varepsilon)}(\theta)=\int_{D^\varepsilon(\theta)} \varphi(v)Dg(\theta+{\eta(\varepsilon)} v)\,dv$$
	$$Df_{\eta(\varepsilon)}(\theta)=\int_{D^\varepsilon(\theta)} \varphi(v)Df(\theta+{\eta(\varepsilon)} v)\,dv.$$
	Therefore, by the definition of Lie bracket we get 
	\begin{multline}\label{bracket1}[f_{\eta(\varepsilon)},g_{\eta(\varepsilon)}](\theta)=\\=\int_{D^\varepsilon(\theta)}\varphi(v)[f,g](\theta)\,dv+E_{f,g}({\eta(\varepsilon)})(\theta)\end{multline}
%	 \begin{multline*}Dg_{\eta(\varepsilon)}(\theta)f_{\eta(\varepsilon)}(\theta)=\\=\int_{D^\varepsilon(\theta)}\varphi(v) Dg(\theta+{\eta(\varepsilon)} v)f_{\eta(\varepsilon)}(\theta)\,dv=\\=\int_{D^\varepsilon(\theta)}\varphi(v)Dg(\theta+{\eta(\varepsilon)} v)f(\theta+{\eta(\varepsilon)} v)\,dv+\\+E_{f,g}({\eta(\varepsilon)})\end{multline*} 
where \begin{multline*}E_{f,g}({\eta(\varepsilon)})(\theta):=\\=\int_{D^\varepsilon(\theta)} \varphi(v)Dg(\theta+{\eta(\varepsilon)} v)[f(\theta+{\eta(\varepsilon)} v)-f_{\eta(\varepsilon)}(\theta)]\,dv-\\-\int_{D^\varepsilon(\theta)} \varphi(v)Df(\theta+{\eta(\varepsilon)} v)[g(\theta+{\eta(\varepsilon)} v)-g_{\eta(\varepsilon)}(\theta)]\,dv.\end{multline*}
If $r>0$, it is easy to show that, if we choose a positive constant $\bar \varepsilon>0$ small enough, we  have  $\theta(q,\sigma,\sqrt{\varepsilon},\tau) \in B_{r}(q)$,  for all $\varepsilon < \bar \varepsilon$. 
 Using the  Lipschitz continuity  of $f$ and $g$ and the convergence of the mollified fields to the original fields, one can easily prove that $|E_{f,g}({\eta})(\theta)|\leq C\eta$,  for  some constant $C$ and  for all $\theta\in B_{r}(q) $.\\  Therefore, by \eqref{bracket1} we get  \begin{multline}
\Big[f_{\eta(\varepsilon)},g_{\eta(\varepsilon)}\Big](\theta(q,\sigma,{\sqrt{\varepsilon}},\tau)) =\\ \int_{D^\varepsilon(\theta(q,\sigma,{\sqrt{\varepsilon}},\tau))}\varphi(v)[f,g](\theta)\,dv \in\\ \overline{\text{co}}\{[f,g](x),\, x \in D \cap B_{r}(q)\}. \end{multline} If we also take $\bar \varepsilon$ small enough so that $\varepsilon < \bar \varepsilon$  implies $\eta=\eta(\varepsilon)<\frac{r}{2C}$, we have proved that for any $r>0$ there is $\bar \varepsilon$ small enough such that \begin{multline*}L(\varepsilon)=\\=\frac{1}{\varepsilon}\int_0^{\sqrt{\varepsilon}}\int_0^{\sqrt{\varepsilon}}\left(\int_{D^\varepsilon(\theta)}\varphi(v)[f,g](\theta)\,dv+E_{f,g}(\eta)\right)d\tau\,d\sigma,\end{multline*}
where we have written $\theta$ instead of $\theta(q,\sigma,{\sqrt{\varepsilon}},\tau)$,
 has distance less than $r$ from $\overline{\text{co}}\{[f,g](x),\, x \in D \cap B_{r}(q)\}$.
\end{proof}
 \addtolength{\textheight}{-12cm}    % This command serves to balance the column lengths

%%%%%%%%%%%%%%%%%%%%%%%%%%%%%%%%%%%%%%%%%%%%%%%%%%%%%%%%%%%%%%%%%%%%%%%%%%%%%%%%

\end{document}